\documentclass[a4paper,oneside,reqno,12pt]{amsart}
\usepackage[utf8]{inputenc}
\usepackage[T1]{fontenc}
\usepackage{amsmath,amssymb,mathtools}
\usepackage{fullpage}
\usepackage[english]{babel}
\usepackage[hidelinks]{hyperref}
\hypersetup{
  pdftitle={A short proof of Ball's plank theorem},
  pdfauthor={Arseniy Akopyan and Alexander Polyanskii}
}

\newtheorem*{coretheorem}{The core theorem}

\DeclareMathOperator{\conv}{conv}
\title{A short proof of Ball's plank theorem}
\author{Arseniy Akopyan}
\author{Alexander Polyanskii}

\address{Arseniy Akopyan, FORA Capital, LLC, Miami, USA}
\email{akopjan@gmail.com}

\address{Alexander Polyanskii, Emory University, 400 Dowman Drive, 30322, Atlanta GA, USA}
\email{apolian@emory.edu}

\begin{document}

\begin{abstract}
We give a short proof of the following theorem of K.~Ball:
\emph{If a centrally symmetric convex body $K\subset \mathbb R^d$ is covered by finitely many planks, then the sum of their relative widths with respect to $K$ is at least~$1$.}
\end{abstract}

\maketitle
\thispagestyle{empty}

A \emph{plank} is the open region between two distinct parallel hyperplanes in $\mathbb R^d$, and its width is the distance between them.
The \emph{width} of a convex body $K\subset \mathbb R^d$ in a given direction is the distance between the two supporting hyperplanes of $K$ perpendicular to that direction. The \emph{relative width} of a plank with respect to $K$ is its width divided by the width of $K$ in the direction perpendicular to the plank's boundary hyperplanes. 

The theorem of K.~Ball~\cite{Ball} states exactly what is said in the abstract. Th.~Bang~\cite{Bang} conjectured this theorem for any convex body, not necessarily centrally symmetric, and this conjecture remains open.
We begin by reformulating the statement in a form that will be more convenient to work with.

\subsection*{Warm-up}

Suppose that a centrally symmetric convex body $K\subset \mathbb R^d$ is covered by $s$ planks with the sum of their relative widths strictly less than $1$. 
Following the original proof, we notice that we can assume that all planks have the same relative width.
We can cover each plank by overlapping parallel planks of relative width $1/n$, increasing their total relative width by at most $s/n$.
Choosing $n$ sufficiently large, we can achieve this covering while keeping the total relative width below $1$. 
By adding extra planks if needed, we can assume that $K$ can be covered by $n$ planks of the same relative width $1/n$.
Next, we reduce the problem to a plank covering problem for a standard crosspolytope in a higher dimension.

Write each of the resulting planks as $\{u:|f_i(u)-c_i|<1\}$, where $f_1,\dots, f_n$ are linear functionals and $c_1,\dots,c_n$ are reals.
Since $K=-K$ and the relative width of the planks is $1/n$, we have $\max_K f_i=n$ for all $i$; choose $u_i\in K$ with $f_i(u_i)=n$. 

For an orthonormal basis $(x_1,\dots,x_n)$ of
$\mathbb R^n$, define a linear map $T:\mathbb R^n\to\mathbb R^d$ by $Tx_i=u_i/n$ and 
the pullback functionals $\ell_i=f_i\circ T$, which clearly satisfy $\ell_i(x_i)=1$.
Notice that
\(
    T(n\conv\{\pm x_1,\dots,\pm x_n\})
  =\conv\{\pm u_1,\dots,\pm u_n\}\subset K.
\), 
The pullback planks $\{x:|\ell_i(x)-c_i|<1\}$ cover the crosspolytope $n\conv\{\pm x_1,\dots,\pm x_n\}$ and in particular the ball $\{x\in \mathbb R^n:\|x\|\leq\sqrt{n}\}$ inscribed in it.
This contradicts the following.

\begin{coretheorem} 
\label{thm:planks}
Let $X:=(x_1,\dots, x_n)$ be an orthonormal basis of $\mathbb R^n$ and let $\ell_1,\dots,\ell_n$ be linear functionals on $\mathbb R^n$ such that $\ell_i(x_i)=1$. Then for any reals $c_1,\dots,c_n$, there is a point $z\in \mathbb R^n$ such that $\|z\|\leq \sqrt{n}$ and $|\ell_i(z)-c_i|\geq1$ for all $i$.
\end{coretheorem}

\begin{proof}   
    Among all orthonormal bases $Y=(y_1,\dots,y_n)$, choose one minimizing
    \[
    F(Y)=\sum\nolimits_i p_i^2,\qquad p_i:=\frac1{\ell_i(y_i)},
    \]
    where $F(Y)=+\infty$ if any denominator vanishes. A minimizer exists by compactness and continuity, and $F(Y)\le F(X)=n$.
    
    For $\varepsilon=(\varepsilon_1,\dots,\varepsilon_n)\in \{\pm 1\}^n$, set 
    $z_\varepsilon=\sum_i\varepsilon_i p_i y_i$. Since $Y$ is orthonormal,
    $\|z_\varepsilon\|^2=\sum_i p_i^2=F(Y)\le n$.
    It remains to show that one of these points $z_\varepsilon$ satisfies $|\ell_i(z_\varepsilon)-c_i|\geq 1$ for all $i$.
    
    For $i\ne j$, rotate $y_i,y_j$ in the plane they span by setting
    $y_i(t)=\cos t\,y_i+\sin t\,y_j$ and
    $y_j(t)=-\sin t\,y_i+\cos t\,y_j$, while fixing the remaining vectors of $Y$. This preserves orthonormality; denote the resulting basis by $Y(t)$. Using the linearity of $\ell_i$ and $\ell_j$, together with the identities
$y_i'(0)=y_j$ and $y_j'(0)=-y_i$, we differentiate $F(Y(t))$ at $t=0$ to obtain
    \begin{equation}
    \label{pi3pj}
    0
    =-2p_i^3\ell_i(y_j)+2p_j^3\ell_j(y_i), \quad \text{so}\quad
    p_i^3\ell_i(y_j)=p_j^3\ell_j(y_i).
    \end{equation}

    If $c_i=0$ for all $i$, choose $\varepsilon$ maximizing
    $Q(\varepsilon)=\sum_i \varepsilon_i p_i^4\ell_i(z_\varepsilon)$.
    If $\varepsilon'$ is obtained by flipping~$\varepsilon_j$, then
    \begin{align*}
    0\le Q(\varepsilon)-Q(\varepsilon')
    &=2\sum\nolimits_{i\ne j}
    \varepsilon_i\varepsilon_j
    \bigl(p_i^4p_j\,\ell_i(y_j) + p_j^4p_i\,\ell_j(y_i) \bigr)\quad\quad 
    \\ &\overset{\eqref{pi3pj}}{=} 4 p_j^4\sum\nolimits_{i\ne j}
    \varepsilon_i \varepsilon_{j} p_i\,\ell_j(y_i)
    = 4p_j^4\bigl(\varepsilon_j\ell_j(z_\varepsilon)-1\bigr),
    \end{align*}
    where the last equality uses $p_j\ell_j(y_j)=1$.
    Therefore $|\ell_j(z_\varepsilon)|\ge1$ for all $j$.

    In the general case, we change the functional to $Q(\varepsilon)=\sum_i \varepsilon_i p_i^4(\ell_i(z_\varepsilon)-2c_i)$. In that case 
    \[
    Q(\varepsilon)-Q(\varepsilon')=4p_j^4\bigl(\varepsilon_j(\ell_j(z_\varepsilon)-c_j)-1\bigr), 
    \]
    and we arrive at $|\ell_j(z_\varepsilon)-c_j|\geq 1$ for all $j$, which finishes the proof of the general case.
\end{proof}

\subsection*{Disclosure of use of AI}
The mathematical ideas of the proof are due to the authors.
ChatGPT was used to assist with the exposition and to simplify some formulas and arguments.

\subsection*{Acknowledgements} AP is supported by the NSF grant DMS 2349045.

\end{document}